\documentclass[11pt]{amsart}
\usepackage{amsmath}
\usepackage{amsfonts}
\usepackage{mathrsfs}
\usepackage{amssymb}
\usepackage{amsthm}
\usepackage[pdftex]{graphicx}
\newtheorem{theorem}{Theorem}[section]
\newtheorem{lemma}[theorem]{Lemma}
\newtheorem{cor}[theorem]{Corollary}
\newtheorem{prop}[theorem]{Proposition}
\theoremstyle{definition}

\theoremstyle{remark}

\numberwithin{equation}{section}
\newcommand{\R}{\mathbb R}

\newcommand{\Z}{\mathbb Z}
\newcommand{\VMO}{\mathrm{VMO}}
\newcommand{\BMO}{\mathrm{BMO}}
\newcommand{\loc}{\mathrm{loc}}
\newcommand{\Sp}{S_{p,p'}}
\newcommand{\avg}{-\hspace{-4mm}\int}
\usepackage[pdftex]{hyperref}
\begin{document}
\author[S.Z. Gautam]{S. Zubin Gautam}
\address{Department of Mathematics, UCLA, Los Angeles, CA 90095-1555, USA.}
\thanks{The author was supported in part by an NSF VIGRE
fellowship.}
\email{sgautam@math.ucla.edu}%
\title{A critical-exponent Balian--Low theorem}
\begin{abstract}
Using a variant of the Sobolev Embedding Theorem, we prove an
uncertainty principle related to Gabor systems that generalizes the
Balian--Low Theorem. Namely, if $f\in H^{p/2}(\R)$ and $\hat f\in
H^{p'/2}(\R)$ with $1<p<\infty$, $\frac{1}{p}+\frac{1}{p'}=1$, then
the Gabor system $\mathcal G(f,1,1)$ is not an exact frame for
$L^2(\R)$.  In the $p=1$ case, we obtain a generalization of the
result in \cite{bcps}.
\end{abstract}
\maketitle
\section{Introduction}\label{intro}
Given a function $f\in L^2(\R)$ and positive constants $\alpha$,
$\beta$, the associated \emph{Gabor system} is \[\mathcal
G(f,\alpha,\beta) := \{e^{2\pi i m\beta \cdot}f(\cdot -
n\alpha)\}_{m,n \in \Z} \subset L^2(\R),\] the collection of
translates and modulates of $f$ by the lattice $\alpha \Z \times
\beta\Z$.  Gabor systems have proven useful in time-frequency
analysis as means for generating orthonormal bases or ``frames" for
$L^2(\R)$. A \emph{frame} for a Hilbert space $\mathcal H$ is a
collection $\{e_n\}\subset \mathcal H$ for which one has the
modified Parseval relation
\begin{equation}A\|x\|_{\mathcal H}^2 \leq \sum_n |\langle x, e_n\rangle|^2
\leq B\|x\|_{\mathcal H}^2 \label{frame} \end{equation} for all
$x\in \mathcal H$ and some \emph{frame constants} $A,B>0$; frames
may be viewed as natural generalizations of orthonormal bases. We
adopt the terminology ``$(A,B)$-frame" for a frame with frame
constants $A$ and $B$.

It is natural to consider under what conditions $\mathcal G(f,
\alpha, \beta)$ generates a frame for $L^2(\R)$; the classical
Balian-Low Theorem is an instance of the uncertainty principle in
this setting (see \textit{e.g.} \cite{daub}):

\begin{theorem}[Balian--Low--Coifman--Semmes]{\label{BLT}}
Let $f\in L^2(\R)$.  If $f\in H^1(\R)$ and $\hat f \in H^1(\R)$,
then $\mathcal G (f,1,1)$ is not a frame for
$L^2(\R)$.\footnote{Known results show that $\alpha=\beta =1$ are
the ``interesting" lattice constants in this setting; see
\textit{e.g.} \cite{daub}.}
\end{theorem}
Here $H^1(\R)$ denotes the usual $L^2$-Sobolev space; thus we see
that if $f$ is suitably well-localized in phase space, then it
cannot generate a Gabor frame. In light of this result, it is
reasonable to ask whether one can alter the regularity assumptions
on $f$ and $\hat f$ to obtain a similar uncertainty principle.

To date, two significant results in this direction have suggested
critical Sobolev regularity assumptions.  The first, essentially due
to Gr\"ochenig \cite{groch}, is:
\begin{theorem}{\label{grochthm}}
Let $1<p,q<\infty$ with $\frac{1}{p} + \frac{1}{q} < 1$.  If $f\in
H^{p/2}(\R)$ and $\hat f\in H^{q/2}(\R)$, then $\mathcal G (f,1,1)$
is not a frame for $L^2(\R)$.
\end{theorem}
From the other direction, Benedetto \textit{et al.} prove the
following in \cite{bcgp}:
\begin{theorem}{\label{bcgpthm}}
Let $\frac{1}{p} + \frac{1}{q} > 1$.  Then there exists a function
$f\in L^2(\R)$ such that that $\mathcal G(f,1,1)$ is a frame (in
fact an orthonormal basis) and such that $f\in H^{p/2}(\R)$ and
$\hat f \in H^{q/2}(\R)$.
\end{theorem}
(In fact, their result is stronger; it allows for stricter
regularity conditions than inclusion in the appropriate Sobolev
spaces.)

Given these results, it is natural to study the critical exponent
case $\frac{1}{p} + \frac{1}{q} = 1$.  In \cite{bcps}, Benedetto
\textit{et al.} conjectured that in fact Theorem \ref{grochthm} can
be extended to this range of exponents, and they proved the
following ``$(1,\infty)$ endpoint" result:

\begin{theorem}{\label{bcpsthm}}
If $f\in H^{1/2}(\R)$ is supported in the interval $[-1,1]$, then
$\mathcal G(f,1,1)$ is not a frame for $L^2(\R)$.
\end{theorem}

The main result of this paper is the following theorem, which
answers the aforementioned conjecture in the affirmative.
\begin{theorem}{\label{mainthm}}
\item
\begin{enumerate}
\item Let $1<p<\infty$.  If $f\in H^{p/2}(\R)$ and $\hat f \in
H^{p'/2}(\R)$, where $\frac{1}{p} + \frac{1}{p'} = 1$, then
$\mathcal G(f,1,1)$ is not a frame for $L^2(\R)$.
\item If $f\in H^{1/2}(\R)$ has compact support, then $\mathcal
G(f,1,1)$ is not a frame for $L^2(\R)$.
\end{enumerate}
\end{theorem}
Note that the $p=2$ case of this theorem is the classical
Balian--Low Theorem (Theorem \ref{BLT}); part 2 is a slight
generalization of Theorem \ref{bcpsthm} (see also Remark (2) below
for a further strengthening of this result).

Before proceeding to the proof of this theorem, we provide some
remarks on its general philosophy in relation to the history of the
problem.  In particular, some discussion of the proof of the
Balian--Low Theorem \ref{BLT} is in order.  The key tool in the
original (incomplete) proof given independently by Balian \cite{bal}
and Low \cite{low} is the \emph{Zak transform}, also known as the
\emph{Weil--Brezin map}. For compactly supported $f\in L^2(\R)$, the
Zak transform $Zf \in L^2_{\mathrm{loc}}(\R^2)$ is given by
\[Zf(x,y) = \sum_{\ell \in \Z}e^{2\pi i \ell y}f(x-\ell).\]
One can view $Zf$ as a function on the unit square $Q_0 := [0,1)
\times [0,1)$, and in fact $Z$ extends to an isomorphism from
$L^2(\R)$ to $L^2(Q_0)$.  We will develop some background on the Zak
transform in Section \ref{background} below. For the present, we
note that to prove either Theorem \ref{BLT} or Theorem
\ref{mainthm}, it suffices to show that
\[\mathrm{ess}\,\inf |Zf| = 0\] under the given regularity
assumptions.\footnote{Here $\mathrm{ess}\, \inf g :=
\inf\big\{\lambda \, \big| \, |\{g\leq\lambda\}|>0\big\}$ is the
essential infimum of $g$, where $|E|$ is the Lebesgue measure of a
set $E$.} Surprisingly, this is the case for \emph{any} function $f$
for which $Zf$ is continuous (see Proposition \ref{zero} below). In
particular, it is worth noting that the proof of this fact is based
on a winding number argument and is hence ``degree-theoretic" in the
topological sense (albeit very simply).

In their original proofs of Theorem \ref{BLT}, Balian and Low
claimed that the regularity conditions $f\in H^1$ and $\hat f\in
H^1$ would force $Zf$ to be continuous; by the remarks above, this
would imply the theorem. However, the regularity conditions only
imply that $Zf\in H^1_{\mathrm{loc}}(\R^2)$, which is not contained
in $C(\R^2)$. This gap in the proof was filled by Coifman and Semmes
and presented in \cite{daub}.

In fact, the Coifman--Semmes argument may be viewed as a simple
prototype of the VMO-degree construction in the Brezis--Nirenberg
theory of \cite{bn1} and \cite{bn2}, which heavily influences our
approach in the current paper.  Broadly speaking, the results of
\cite{bn1} and \cite{bn2} show that in many cases VMO maps are as
good as continuous maps for the purposes of degree theory. (Here
$\VMO(\R^n)$ is Sarason's space of functions of \emph{vanishing mean
oscillation} on $\R^n$; see Section \ref{background} below.) In
accordance with this principle, Coifman and Semmes first prove that
under the given regularity assumptions $Zf \in \VMO (\R^2)$; their
argument gives the $n=p=2$ case of the following endpoint Sobolev
embedding theorem (see \textit{e.g.} \cite{bn1}):
\begin{theorem}{\label{sobolev}}
Let $1\leq p < \infty$, and let $s = p/n$.  Then $W^{s,p}(\R^n)
\subset \VMO(\R^n)$ with continuous embedding, where $W^{s,p}$ is
the usual $L^p$-Sobolev space.
\end{theorem}
This fact is then used to run a modified winding number argument and
prove that $\mathrm{ess}\,\inf |Zf| = 0$, from which the theorem
follows. (The Coifman--Semmes proof as presented in \cite{daub} does
not explicitly mention VMO, BMO or the above Sobolev embedding, but
the methods are present without the terminology.)

For the proof of our main result, Theorem \ref{mainthm}, we take a
parallel approach.  As noted, prior to the results of \cite{bcps},
the best known result was Theorem \ref{grochthm}.  This latter
follows from the results of \cite{groch}, in which it is shown that
under the given regularity assumptions $f$ belongs to the Wiener
algebra
\[W(\R) = \{f\in C(\R)\, | \, \sum_{k\in \Z} \sup_{x\in [0,1]}
|f(x+k)| < \infty\}.\] This in turn immediately implies that the Zak
transform of $f$ is continuous. However, the proof is invalid for
the critical regularity case where $\frac{1}{p}+\frac{1}{q} = 1$, so
we expect that $Zf$ ``barely" fails to be continuous under our
regularity assumptions.  Thus it seems reasonable to expect that
$Zf\in \VMO(\R^2)$; that this is in fact true is the key step of our
proof, established by a variant of the above Sobolev embedding
theorem (Theorem \ref{sobolevvariant} below). We combine this with a
simplified version of the Coifman--Semmes winding number argument
(essentially drawn from \cite{bn2}) to yield the final result.

In the sequel, we will write ``$A \lesssim B$" if $A \leq c B$ for
some universal constant $c$; ``$A \sim B$" means $A \lesssim B
\lesssim A$.  Subscripts on the symbols ``$\lesssim$" and ``$\sim$"
will denote dependence of the implied constants.

\section{Background and preliminaries:  The Zak transform and VMO}{\label{background}}
We begin by recalling some basic facts about the Zak transform.  As
stated above, for compactly supported $f\in L^2(\R)$, the Zak
transform of $f$ is defined (almost everywhere) by \[Zf(x,y) =
\sum_{\ell \in \Z} e^{2\pi i \ell y} f(x-\ell).\]  It is easily seen
that $Zf$ verifies the ``quasi-periodicity" relations
\begin{eqnarray}
Zf(x+1,y) & = & e^{2\pi i y} Zf(x,y) \nonumber \\
Zf(x,y+1) & = & Zf(x,y), \label{qp}
\end{eqnarray}
so that $Zf$ is completely determined by its values on the unit cube
$Q_0 \subset \R^2$.  As mentioned above, $Z$ actually extends to a
unitary isomorphism from $L^2(\R)$ to $L^2(Q_0)\cong L^2(\mathbb
T^2)$.  This can easily be seen by examining its action on the
orthonormal basis $\{e_{m,n}\}$ of $L^2(\R)$, where \[e_{m,n}(x) =
e^{2\pi i n x} \;\mathbf 1_{[0,1)} (x - m),\;\;m,n\in\Z;\] this
basis is mapped to the usual Fourier basis of $L^2(\mathbb T^2)$ by
the Zak transform.  Thus we may view $Z$ as a map from $L^2(\R)$ to
either $L^2(\mathbb T^2)$ or $L^2_{\loc}(\R^2)$.

$Zf$ provides a time-frequency representation of $f$; in fact,
viewed as an element of $L^2(\mathbb T^2)$, $Zf$ is the Fourier
transform of the \emph{Gabor coefficients} $(f_{m,n}) \in
\ell^2(\Z\times\Z)$, defined by \[f_{m,n} = \langle f, e_{m,n}
\rangle.\]  Similarly, the Zak transform is intimately connected
with the frame properties of Gabor systems.
\begin{prop}{\label{AB}}
Let $f\in L^2(\R)$.  Then $\mathcal G(f,1,1)$ is an $(A,B)$-frame
for $L^2(\R)$ if and only if $A^{1/2} \leq |Zf| \leq B^{1/2}$ almost
everywhere.
\end{prop}
This is complemented by the following somewhat curious fact, as
mentioned above.
\begin{prop}{\label{zero}}
If $f\in L^2(\R)$ has continuous Zak transform, then $Zf$ must have
a zero.
\end{prop}
For the proofs of these results and more on the Zak transform, see
\textit{e.g.} \cite{daub} and \cite{fol}.  In light of Proposition
\ref{AB}, we see that in order to prove an obstruction result such
as the Balian--Low Theorem \ref{BLT} or Theorem \ref{mainthm}, it
suffices to show that $\mathrm{ess}\,\inf |Zf| = 0$.  We will
accomplish this in part by proving an analogue of Proposition
\ref{zero} (Proposition \ref{vmozero} below).

We now discuss the regularity properties of the Zak transform of a
function $f$ satisfying some given time-frequency localization (or
regularity) conditions.  For convenience, we introduce the notation
$S_{p,q}$ with $0<p,q<\infty$ for the Hilbert space
\[S_{p,q}:= \big\{g\in L^2(\R^2)\: \big |\: \int_{\R^2} |\hat g(\xi_1,\xi_2)|^2
(1+|\xi_1|^p + |\xi_2|^{q}) \, d\xi_1 \,d\xi_2 < \infty\big\},\]
equipped with the norm \[\|g\|_{S_{p,q}} = \bigg( \int_{\R^2} |\hat
g(\xi_1, \xi_2)|^2 (1+|\xi_1|^p + |\xi_2|^{q}) \, d\xi_1 \,d\xi_2 \,
\bigg)^{1/2}.\] $S_{p,q}$ should be thought of as a modified Sobolev
space; when $p=q$, $S_{p,p}$ coincides with the usual inhomogeneous
Sobolev space $H^{p/2}(\R^2)$, with equivalent norms.

The Zak transform of a function $f$, being a time-frequency
representation of $f$, naturally inherits the smoothness properties
of $f$ and $\hat f$ in the following sense.
\begin{lemma}{\label{regularity}}
Let $f\in H^{s_1}(\R)$ and $\hat f \in H^{s_2}(\R)$ with
$s_1,s_2>0$.  Then for any smooth, compactly supported function
$\varphi \in C_c^\infty(\R^2)$, we have $\varphi Zf \in
S_{2s_1,2s_2}$.
\end{lemma}
\noindent \textbf{Proof}  For $j=1,2$, we write $\nabla_j$ for the
$j$-th distributional partial derivative operator on the space of
tempered distributions $\mathcal S'(\R^2)$; $\nabla$ denotes the
distributional derivative on $\mathcal S'(\R)$.  Similarly, for
$s\geq 0$, we define the inhomogeneous fractional derivatives
$\langle \nabla_j \rangle^s$ as Fourier multipliers  on $\mathcal
S'(\R^2)$ with symbols $\langle \xi_j\rangle^s$.\footnote{Here
$\langle a \rangle = (1+|a|^2)^{1/2}$.}

Recalling the definition \[Zf(x_1,x_2) = \sum_{\ell \in \Z} e^{2\pi
i \ell x_2}f(x_1 - \ell)\] for compactly supported $f$, it is easy
to check that for all $f\in H^{s}(\R)$ \[Z ( \nabla^n f ) =
\nabla_1^n (Zf)\] in the sense of distributions for any nonnegative
integer $n\leq s$.

Now for $j =1,2$ and $s\geq 0$, let $H_j^s(\R^2)$ denote the
modified Sobolev space
\[H_j^s(\R^2) = \{f\in L^2(\R^2) \;|\; \|f\|_{H_j^s} = \|\langle
\nabla_j\rangle^s f\|_2 <\infty\}.\]  Fix a compactly supported bump
function $\varphi \in C_c^\infty(\R^2)$.  When $k$ is an integer,
the Leibniz rule for weak derivatives yields
\begin{eqnarray*}\|\varphi Zf\|_{H^k_1} & \lesssim_{k}\phantom{_{,\varphi}} &
\|\varphi Zf\|_{L^2(\R^2)}\, + \,\sum_{m=0}^k \|( \nabla_1^m
\varphi) ( \nabla_1^{k-m} Zf) \|_2
\\ & =\phantom{_{k,\varphi}} & \|\varphi Zf\|_{L^2(\R^2)} \, + \, \sum_{m=0}^k
\|(\nabla_1^m \varphi )Z(\nabla^{k-m} f) \|_2 \\
& \lesssim_{k,\varphi} & \sum_{m=0}^k \|Z(\nabla^m f)
\|_{L^2(\mathrm{supp}(\varphi))} \\ &
\lesssim_{\varphi}\phantom{_{,k}}& \|f\|_{H^k(\R)}.
\end{eqnarray*}  The last two inequalities follow from the quasi-periodicity
relations (\ref{qp}) and the unitarity of $Z$ viewed as a map into
$L^2(\mathbb T^2)$, which imply that \[\|Zg\|_{L^2(K)} \lesssim_K
\|g\|_{L^2(\R)}\] for all compact $K\subset \R^2$.  Thus $\varphi Z$
is a bounded linear operator from $H^k(\R)$ to $H^k_1(\R^2)$ for
integer values of $k$.
 The spaces
$H^s_1(\R^2)$ can be interpolated via the complex method as with the
traditional Sobolev spaces (see \textit{e.g.} \cite{af}), so we
obtain in fact that $\varphi Z$ is bounded from $H^s(\R)$ to
$H^s_1(\R^2)$ for all $s\geq 0$. (Equivalently, one can work on the
Fourier transform side and appeal to Stein's weighted interpolation
theorem; see \cite{stein2}.)

A similar argument also shows that \[ \langle \nabla_2 \rangle^s
(\varphi Zf) \in L^2(\R^2)\] whenever $\hat f \in H^s(\R)$, once one
applies the well-known fact that
\[Z\hat f(x,y) = e^{2\pi i xy}Zf(-y,x),\] which follows from the
Poisson Summation Formula. So when $f\in H^{s_1}(\R)$ and $\hat f
\in H^{s_2}(\R)$, we have \[\big( \langle \nabla_1 \rangle^{s_1} +
\langle \nabla_2 \rangle^{s_2} \big) (\varphi Zf) \in L^2(\R^2)\]
for all $\varphi \in C_c^\infty(\R^2)$.  The lemma then follows
immediately from Plancherel's Theorem.

\newpage
Finally, we recall some basic facts about the space $\VMO(\R^n)$.
Recall that $\BMO(\R^n)$ is the space of functions (modulo
constants) of bounded mean oscillation on $\R^n$,
\[\BMO(\R^n) = \big\{ f\in L^1_\loc (\R^n) \; \big| \; \|f\|_\BMO :=
\sup_Q \Big( - \hspace{-3.8mm}{\int}_Q |f(x) - -
\hspace{-3.4mm}{\int}_Q f| \, dx\Big) \, < \, \infty\, \big\},\]
where the supremum is taken over cubes $Q$ in $\R^n$, and
\[ \avg_E g := \frac{1}{|E|} \int_E g\] denotes the average of a
function $g$ over a Lebesgue measurable set $E$. We define
$\VMO(\R^n)$ to be the closure of the uniformly continuous functions
in the BMO-norm.  We will also use a more concrete characterization
of VMO:  $f\in \BMO(\R^n)$ is in VMO if and only if
\[ \lim_{a\rightarrow 0} \, \sup_{|Q| \leq a} \Big( -\hspace{-4.2mm}{\int}_Q |f(x) -
\, \avg_Q f|\, dx \Big) = 0.\]  For the proof of this and other
equivalent characterizations of VMO, see \cite{sar}.

\section{The embedding into VMO}
As a first step, we establish that whenever $f\in H^{p/2}(\R)$ and
$\hat f \in H^{p'/2}(\R)$ for $1<p,p'<\infty$ with $p$ and $p'$
conjugate, we have $Zf\in \VMO(\R^2)$.  The key step of our proof is
the following analogue of the endpoint Sobolev embedding in Theorem
\ref{sobolev} for the spaces $S_{p,q}$.

\begin{theorem}{\label{sobolevvariant}}
Let $1<p<\infty$, and let $p'$ be the conjugate exponent to $p$.
Then $S_{p,p'} \subset \VMO(\R^2)$ with \[ \|f\|_{\BMO} \lesssim
\|f\|_{\Sp}.\]
\end{theorem}
\noindent \textbf{Proof}  Without loss of generality, we assume
$1<p<2$, so that $p'>2>p$.  (The case $p=2$ is a special case of
Theorem \ref{sobolev}, as mentioned above.)  We will use the
Littlewood--Paley characterization of BMO.

Let $\{\psi_k\}_{k\in \Z}$ be a Littlewood--Paley partition of unity
on the frequency space $\R^2$, so that each $\psi_k$ is a
nonnegative smooth bump function supported on the annulus $\{2^{k-1}
\leq |\xi| \leq 2^{k+1}\}$, with \[\sum_k \psi_k(\xi) = 1, \: a.e.
\, \, \xi\in\R^2.\]  Let $P_k$ denote the corresponding
Littlewood--Paley projection operators, so that each $P_k$ is a
Fourier multiplier with symbol $\psi_k$.  Then for all $f\in
L^2(\R^2)$, we have \[\|f\|_{\BMO(\R^2)} \sim_c \sup_Q
\bigg(-\hspace{-4.2mm}{\int}_Q\, \sum_{k\geq -\log_2\ell(Q) + c}
|P_kf|^2 \bigg)^{1/2},\] where the supremum is taken over cubes
$Q\subset \R^2$ of dyadic side lengths $\ell(Q)$, and $c\geq 0$. For
our purposes, it suffices to take $c=3$.  (This is essentially a
discrete version of Theorem 3 in Chapter IV, \S 4.3 of \cite{stein};
see also \S 4.5 of the same.)

Let $f\in L^2(\R^2)$, and fix a cube $Q \subset \R^2$ with $\ell(Q)
= 2^{-k_0}$. As a first step, we prove the estimate
\begin{equation}{\label{piece}}
-\hspace{-4.2mm}{\int}_Q |P_kf|^2 \lesssim \int_{\R^2} |\widehat{P_k
f}(\xi)|^2 (1 + |\xi_1|^p + |\xi_2|^{p'})\, d\xi_1\, d\xi_2
\end{equation}
for each Littlewood--Paley piece $P_kf$ with $k \geq k_0 + 3$. Let
$\psi \in \mathcal S(\R^2)$ be a nonnegative Schwartz function
adapted to $Q$, such that $\psi \geq 1$ on $Q$ and $\hat \psi$ is
supported on the cube of length $2^{k_0}= \ell(Q)^{-1}$ centered at
$0$.  Then $\psi$ will satisfy the estimate \[\|\psi\|_2 \lesssim
2^{-k_0},\] with implied constant independent of $Q$. Then we have
\begin{equation}{\label{cubedecomp}} -\hspace{-4.1mm}{\int}_Q |P_kf|^2 \leq 2^{2k_0}\|\psi
P_kf\|_{L^2(\R^2)}^2 = 2^{2k_0}\sum_{\substack{\textrm{$J$ dyadic}
\\ \ell(J) = 2^{k_0}}} \|\hat \psi
* \widehat{P_kf}\|_{L^2(J)}^2,
\end{equation}
where the last sum is taken over the disjoint cubes $J$ in the
dyadic mesh at scale $2^{k_0}$.  Let $w_p$ be the weight function
\[w_p(\xi) = (1+|\xi_1|^p + |\xi_2|^{p'})^{1/2}.\]  Then for each $J$,
Young's inequality and Cauchy--Schwarz yield \[ \|\hat \psi
* \widehat{P_k f}\|_{L^2(J)} \leq \|\hat \psi\|_2\,
\|\widehat{P_kf}\|_{L^1(3J)} \lesssim 2^{-k_0} \|\widehat{P_kf}
\cdot w_p\|_{L^2(3J)} \, \|w_p^{-1}\|_{L^2(3J)}\] since $\hat \psi$
is supported on a cube of side length $2^{k_0} = \ell (J)$ (here
$3J$ denotes the cube with the same center as $J$ and $\ell(3J) =
3\ell(J)$).  Summing over $J$ as in (\ref{cubedecomp}) yields \[
\avg_Q |P_kf|^2 \lesssim \sum_{J\in\, \mathcal J_{k_0}}
\|\widehat{P_kf}\cdot w_p \|_{L^2(3J)}^2 \|w_p^{-1}\|_{L^2(3J)}^2.\]
Here we write $\mathcal J_{k_0}$ for the collection of ``admissible"
cubes on which $\widehat{P_k f}$ does not vanish identically for
every $k \geq k_0+3$. In particular, since $\widehat{P_kf}$ is
supported on the annulus $\{2^{k-1} \leq |\xi| \leq 2^{k+1}\}$,
$\mathcal J_{k_0}$ omits cubes sufficiently close to the origin
(\textit{viz.}, the shaded cubes in Figure \ref{pic} below).

\begin{figure}
\setlength{\unitlength}{.9 cm}
\begin{picture}(12,12)(-6,-6)
\put(5.2,-.1){{$\xi_1$}} \put(-.1,5.2){{$\xi_2$}}
\put(1,0){\line(0,1){5}}\put(1,0){\line(0,-1){5}}
\put(2,0){\line(0,1){5}}\put(2,0){\line(0,-1){5}}
\put(3,0){\line(0,1){5}}\put(3,0){\line(0,-1){5}}
\put(4,-5){\line(0,1){5}}\put(-4,-5){\line(0,1){10}}
\put(4,0){\line(0,1){1.3}}\put(4,2.2){\line(0,1){2.8}}
\put(-1,0){\line(0,1){5}}\put(-1,0){\line(0,-1){5}}
\put(-2,0){\line(0,1){5}}\put(-2,0){\line(0,-1){5}}
\put(-3,0){\line(0,1){5}}\put(-3,0){\line(0,-1){5}}
\put(-5,1){\line(1,0){10}}\put(-5,2){\line(1,0){8.5}}
\put(-5,3){\line(1,0){10}}\put(-5,-1){\line(1,0){10}}
\put(-5,-2){\line(1,0){10}}\put(-5,-3){\line(1,0){10}}
\put(-5,4){\line(1,0){10}}\put(-5,-4){\line(1,0){10}}
\put(2.35,.35){{$J^*$}}
\put(2.1,-.4){{\scriptsize{$2^{k_0+1}$}}}%
\qbezier(0,0)(2,1.2)(4.5,1.5)%
\put(3.8,1.7){{\footnotesize{$|\xi_1|^p = |\xi_2|^{p'}$}}}%
\setlength{\unitlength}{.45 cm}%
\multiput(-4,3)(1,0){8}{\line(1,1){1}}%
\multiput(-4,2)(1,0){8}{\line(1,1){1}}%
\multiput(-4,1)(1,0){8}{\line(1,1){1}}%
\multiput(-4,0)(1,0){8}{\line(1,1){1}}%
\multiput(-4,-1)(1,0){8}{\line(1,1){1}}%
\multiput(-4,-2)(1,0){8}{\line(1,1){1}}%
\multiput(-4,-3)(1,0){8}{\line(1,1){1}}%
\multiput(-4,-4)(1,0){8}{\line(1,1){1}}%

\setlength{\unitlength}{.9 cm}
\linethickness{.4 mm}%
\put(-2,1){\line(1,0){4}} \put(-2,2){\line(1,0){4}}
\put(-2,-1){\line(1,0){4}} \put(-2,-2){\line(1,0){4}}
\put(-2,-2){\line(0,1){4}} \put(-1,-2){\line(0,1){4}}
\put(1,-2){\line(0,1){4}} \put(2,-2){\line(0,1){4}}
\put(0,0){\vector(1,0){5}} \put(0,0){\vector(0,1){5}}
\put(0,0){\vector(-1,0){5}} \put(0,0){\vector(0,-1){5}}
\end{picture}
\caption{Inadmissible cubes for $\mathcal J_{k_0}$, cube $J^*$ (at
scale $k_0 \gg 1$). \label{pic}}
\end{figure}
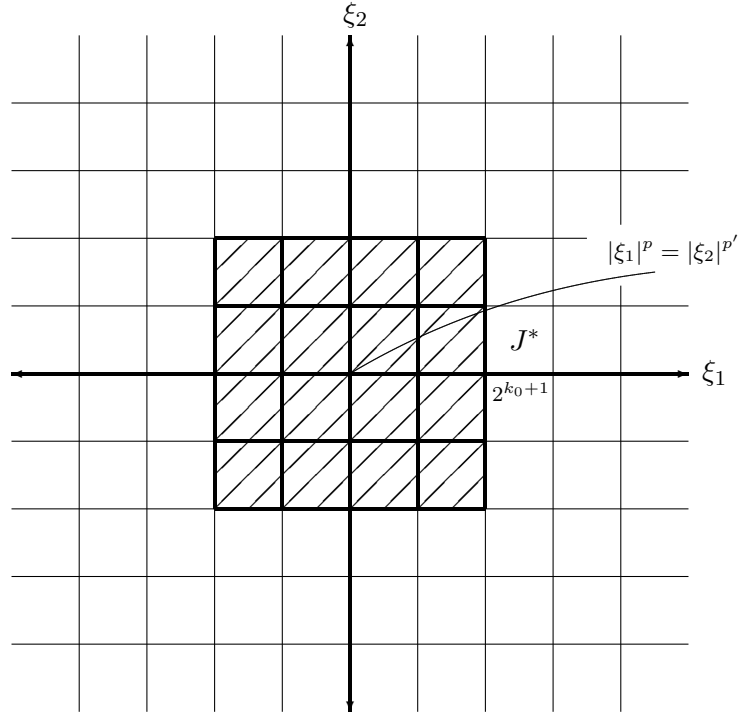

It so happens that this is enough to ensure that for all $J\in
\mathcal J_{k_0}$ \begin{equation}{\label{weight}}
\|w_p^{-1}\|_{L^2(3J)} \lesssim 1,\end{equation} where the implied
constant is \emph{independent} of the scale parameter $k_0$. By
symmetry and monotonicity considerations on $w_p$, and since $p<p'$,
it suffices to bound the term corresponding to $J^* = [2\cdot
2^{k_0}, 3\cdot 2^{k_0}) \times [0, 2^{k_0})$ (see Figure
\ref{pic}).  Now we have the estimate
\begin{eqnarray*}
\|w_p^{-1}\|_{L^2(3J^*)}^2 & = & \int_{2^{k_0}}^{2^{k_0+2}}
\int_{-2^{k_0}}^{2^{k_0+1}} \frac{d\xi_2\, d\xi_1}{1+|\xi_1|^p +
|\xi_2|^{p'}}\\
& \lesssim & \int_{2^{k_0}}^{2^{k_0+2}} \int_0^{2^{k_0+1}}
\frac{d\xi_2\, d\xi_1}{1+|\xi_1|^p + |\xi_2|^{p'}}\\
& \leq & \int_{2^{k_0}}^{2^{k_0+2}} \int_0^{\xi_1^{p/p'}} \xi_1^{-p}
d\xi_2\, d\xi_1 \; + \; \int_{2^{k_0}}^{2^{k_0+2}}
\int_{\xi_1^{p/p'}}^{2^{k_0+1}}\xi_2^{-p'} d\xi_2\, d\xi_1\\
& \lesssim & 1 + 2^{k_0 (2-p')},
\end{eqnarray*}
which is bounded for $k_0\geq 0$ as $p'>2$.  For $k_0 < 0$, we
simply note that $|3J^*| \lesssim 1$ and $\|w_p^{-1}\|_\infty = 1$,
and we obtain (\ref{weight}).  This implies \[\avg_Q |P_kf|^2
\lesssim \sum_{\substack{\textrm{$J$ dyadic}
\\ \ell(J) = 2^{k_0}}} \|\widehat{P_kf} \cdot w_p\,
\|_{L^2(3J)}^2 \lesssim \|\widehat{P_kf}\cdot w_p\,
\|_{L^2(\R^2)}^2,\] which is the desired estimate (\ref{piece}) on
$P_kf$.  Summing this estimate in $k\geq k_0 +3$ and taking the
supremum over all cubes $Q$, we obtain the BMO estimate \[\|f\|_\BMO
\lesssim \|f\|_{\Sp}.\]  But since Schwartz functions are dense in
$\Sp$, we actually have $\Sp \subset \VMO(\R^2)$, as VMO is the
BMO-closure of the uniformly continuous functions.  This concludes
the proof of the theorem.

\bigskip
From this and Lemma \ref{regularity} (combined with the
quasi-periodicity property (\ref{qp})) we obtain:
\begin{cor}{\label{p>1}}
For $1<p<\infty$, if $f\in H^{p/2}(\R)$ and $\hat f \in
H^{p'/2}(\R)$, then $Zf\in \VMO(\R^2)$.
\end{cor}

We now turn to the ``endpoint" regularity case where $p=1$; of
course, the dual localization condition ``$\hat f \in H^{\infty
/2}$" requires suitable interpretation.  For our purposes, as in
Theorem \ref{mainthm}.2, we will take this to mean $f$ has compact
support; this is less restrictive than the condition
$\mathrm{supp}\, (f) \subset [-1,1]$ in Theorem \ref{bcpsthm}. In
this setting, we will show directly that $Zf\in \VMO(\R^2)$,
provided that $Zf\in L^{\infty}(\R^2)$. (The additional boundedness
assumption on $Zf$ will be acceptable for our purposes, in light of
Proposition \ref{AB}.)

\begin{lemma}{\label{p=1}}
Suppose $f\in H^{1/2}(\R)$ has compact support, with $Zf\in
L^\infty(\R^2)$.  Then $Zf\in\VMO(\R^2)$.
\end{lemma}
\noindent \textbf{Proof}  Fix a large cube $Q^* \subset \R^2$ such
that the unit cube $Q_0$ is contained in the interior of $Q^*$.  By
the quasi-periodicity relations ($\ref{qp}$) for the Zak transform,
it suffices to prove $Zf\in \VMO(Q^*)$, in the sense that
\[\lim_{a\rightarrow 0} \; \sup_{\substack{|Q|< a, \\ Q\subset Q^*}}
\Big( -\hspace{-4.15mm}{\int}_Q| Zf(x) - \avg_Q Zf|\, dx \Big) =
0.\] But since $f$ has compact support, its Zak transform is a
\emph{finite} sum
\[Zf(x,y) = \sum_{|\ell| \leq N} e^{2\pi i \ell y}f(x-\ell)\] for
$(x,y) \in Q^*$.  We have $\|f\|_\infty \leq \|Zf\|_\infty<\infty$;
this can be seen for instance by fixing $x$ and viewing $Zf(x,y)$ as
a Fourier series in $y$.  Since $f\in H^{1/2}(\R)$, we also have
$f\in \VMO(\R)$ by Theorem \ref{sobolev}.

A simple calculation shows that if $g,h \in \VMO(\R) \cap L^\infty
(\R)$, then their tensor product $g \otimes h$ lies in $\VMO(\R^2)$.
The restriction of $Zf$ to $Q^*$ agrees with a finite sum of such
tensor products, so we have $Zf\in \VMO(Q^*)$, and the lemma
follows.

\section{The winding number argument}
Recall that the Zak transform of a function $f$ satisfies the
quasi-periodicity relations (\ref{qp}):
\begin{eqnarray*}
Zf(x+1,y) & = & e^{2\pi i y} Zf(x,y)\\
Zf(x,y+1) & = & Zf(x,y)
\end{eqnarray*}
for a.e. $(x,y) \in \R^2$; as mentioned before, any continuous
function satisfying these relations must have a zero.  In fact, the
same is essentially true of bounded VMO functions.

\begin{lemma}{\label{vmozero}}
Suppose $F\in \VMO(\R^2) \cap L^\infty(\R^2)$ satisfies
\begin{eqnarray*}
F(x+1,y) & = & e^{2\pi i y}F(x,y)\\
F(x,y+1) & = & F(x,y)
\end{eqnarray*}
almost everywhere.  Then $\mathrm{ess}\, \inf |F| = 0$.
\end{lemma}
\noindent \textbf{Proof}\footnote{This argument is almost identical
to that of Coifman and Semmes in \cite{daub}; the difference is that
they use a \emph{quantitative} VMO estimate, whereas we need to make
do with only the qualitative assumption that $F\in \VMO$.  We
present the whole argument for self-containment.} We proceed by
contradiction. By scaling, we may assume that there exists $d>0$
such that \[d \leq |F| \leq 1\] almost everywhere. Let
$Q_\varepsilon (x,y)$ denote the cube of side length $\varepsilon$
centered at $(x,y)\in \R^2$, and define
\[F_\varepsilon (x,y) := \avg_{Q_\varepsilon(x,y)}
F,\] the average of $F$ over $Q_\varepsilon (x,y)$. $F_\varepsilon$
is continuous and satisfies the modified quasi-periodicity relations
\begin{eqnarray} F_\varepsilon(x+1,y) & = & e^{2\pi i y}
F_\varepsilon(x,y) + \Phi_\varepsilon(x,y) \nonumber\\
F_\varepsilon(x,y+1) & = & F_\varepsilon(x,y),
\label{modqp}\end{eqnarray} where the error term $\Phi_\varepsilon$
satisfies
\[|\Phi_\varepsilon (x,y)| \lesssim \varepsilon.\]  Moreover, for
$\varepsilon$ sufficiently small, $F_\varepsilon$ is also bounded
from below, as we now show.  Since $F\in \VMO(\R^2)$, we may choose
$\varepsilon_0$ so that \[ \avg_{Q_\varepsilon (x,y)} |F -
F_\varepsilon(x,y)| \leq \frac{d}{2}\] for all $(x,y) \in \R^2$ and
$\varepsilon < \varepsilon_0$. Then simply applying the triangle
inequality, we have
\begin{eqnarray*}
|F_\varepsilon(x,y)| & = & \avg_{Q_\varepsilon (x,y)}
|F_\varepsilon(x,y)|\\
& \geq & \avg_{Q_\varepsilon(x,y)}|F| - \avg_{Q_\varepsilon(x,y)} |F
- F_\varepsilon(x,y)| \; \geq \; \frac{d}{2},
\end{eqnarray*}
since we assume $|F|\geq d$ almost
everywhere. Thus $\frac{d}{2} \leq |F_\varepsilon| \leq 1$ almost
everywhere for $\varepsilon < \varepsilon_0$.

However, this is impossible, as the relations (\ref{modqp}) force
the curve $\Gamma_\varepsilon := F_\varepsilon (\partial Q_0)$ to
have nonzero winding number about $0$ for $\varepsilon$ sufficiently
small, where $Q_0 = [0,1] \times [0,1]$ is the unit cube in $\R^2$.
To make this contradiction more precise, we give the same argument
as Coifman and Semmes.  Note that since $F_\varepsilon$ is
continuous with $\frac{d}{2} \leq |F_\varepsilon| \leq 1$, we can
define a continuous branch $\gamma_\varepsilon$ of $\log
F_\varepsilon$. From the modified quasi-periodicity conditions
(\ref{modqp}), we have
\begin{eqnarray}\gamma_\varepsilon (x+1, y) & = & \gamma_\varepsilon
(x,y) + 2\pi i j + 2\pi i y + \Psi_\varepsilon (x,y)\nonumber \\
\gamma_\varepsilon (x, y+1) & = & \gamma_\varepsilon (x,y) + 2\pi i
k\end{eqnarray} for all $x,y$ in some simply connected neighborhood
$U$ of $Q_0$.  Here $j,k\in \Z$ are constant on $U$ by continuity of
$\gamma_\varepsilon$, and
\[|\Psi_\varepsilon| \leq - \log \bigg(1-
\frac{|\Phi_\varepsilon|}{|F_\varepsilon|}\bigg) \lesssim
\frac{|\Phi_\varepsilon|}{|F_\varepsilon|}\] provided that
$|\Phi_\varepsilon| / |F_\varepsilon|$ is sufficiently small.  This
can be arranged by taking $\varepsilon$ sufficiently small, since
$|\Phi_\varepsilon| \lesssim \varepsilon$ and $|F_\varepsilon| \geq
d/2$; thus for $\varepsilon$ small we have \[|\Psi_\varepsilon| <
1\] on $U$.  To obtain the contradiction, we simply compute
\begin{eqnarray*}
0 & = & \phantom{+\;} \gamma_\varepsilon (1,0) - \gamma_\varepsilon
(0,0) + \gamma_\varepsilon (1,1) - \gamma_\varepsilon(1,0)\\ & & +\;
\gamma_\varepsilon (0,1) - \gamma_\varepsilon(1,1) +
\gamma_\varepsilon (0,0) - \gamma_\varepsilon(0,1)\\
& = & \phantom{+\,}\Psi_\varepsilon (0,0) - \Psi_\varepsilon (0,1) -
2\pi i \, \neq \, 0,
\end{eqnarray*}
since $|\Psi_\varepsilon|<1$.  Thus our original assumption that
$|F| \geq d$ almost everywhere must be false, and hence
$\mathrm{ess} \, \inf |F| = 0$ as desired.

\bigskip
From this, we can deduce our main result, Theorem \ref{mainthm}.
Suppose $f\in L^2(\R)$ satisfies either of the prescribed
time-frequency regularity conditions, and suppose furthermore that
$\mathcal G(f,1,1)$ is an $(A,B)$-frame for $L^2(\R)$.  Then by
Proposition \ref{AB}, we have \[A^{1/2} \leq |Zf| \leq B^{1/2} \;\;
a.e.\]  Moreover, by either Corollary \ref{p>1} or Lemma \ref{p=1},
$Zf\in \VMO(\R^2) \cap L^\infty(\R^2)$.  But, by Lemma
\ref{vmozero}, this is impossible; this contradiction concludes the
proof of the theorem.

\section{Remarks and acknowledgments}
\begin{enumerate}
\item The second part of the proof of Lemma \ref{vmozero}
essentially shows that the continuous maps $F_\varepsilon$ all have
nonzero degree at the point $0$, for $\varepsilon$ sufficiently
small.  In the context of \cite{bn1} and \cite{bn2}, this is a
manifestation of the stability of degree under VMO-convergence.  In
fact, the (integer-valued) VMO degree of $F$ at a point $p$ is
\emph{defined} as
\[\VMO \textrm{-deg } (F,p) := \mathrm{deg}(F_\varepsilon, p) \; \;
\textrm{for $\varepsilon < \varepsilon_0$},\] up to some domain
considerations.  As mentioned before, the $H^1(\R^2)$ argument of
Coifman and Semmes can be viewed as a prototype of the
Brezis--Nirenberg theory in a relatively simple case; for more on
the VMO degree theory and related topics, see \textit{e.g.}
\cite{brezis}, \cite{bn1}, \cite{bn2}, and \cite{bbm}.
\item The localization condition that $f$ be compactly supported in
Theorem \ref{mainthm}.2 is far from sharp.  An inspection of the
proof of Lemma \ref{p=1} shows that it is sufficient to demand
$\BMO$-convergence of the series defining the Zak transform near the
unit square, which would in turn be implied by uniform convergence
of the series.  This latter can easily be guaranteed by simply
requiring the mild decay condition
\[\sum_{\ell\in \Z}\|f\chi_{[\ell, \ell+1]}\|_{L^\infty} <\infty,\] for
instance.  This observation is due to the author, C. Heil, and A.
Powell, arising from a joint discussion.  However, as noted, this
decay condition is stronger than is necessary to guarantee
$\operatorname{Z} f \in \VMO$.
\item Theorem \ref{sobolevvariant} is of mild
interest in its own right.  The spaces $S_{p,q}$ above were chosen
\textit{ad hoc} for the setting of the Balian-Low Theorem; one might
hope to prove an embedding result for spaces with analogous
$L^r$-based regularity conditions, $r\neq 2$.
\end{enumerate}

\bigskip
The author thanks Terence Tao, Rowan Killip, and especially
Christoph Thiele for much valuable advice.  Thanks also go to Chris
Heil and Alex Powell for the discussion mentioned in Remark (2), as
well as to Victor Lie for a useful discussion on the same point.

\end{document}